\documentclass{elsart}
\usepackage{amssymb}
\usepackage{amsmath}
\usepackage{epsfig}
\usepackage{amsfonts}
\usepackage{latexsym}

\begin{document}
\begin{frontmatter}
\title{Binary Matrices under the Microscope: \\
A Tomographical Problem}

%
%\author{Andrea Frosini\inst{1} \and Maurice Nivat\inst{2}}
%\institute{Dipartimento di Scienze Matematiche ed Informatiche
%``Roberto Magari" Universit\`a degli Studi di Siena, Pian dei
%Mantellini 44, 53100, Siena, Italy \\
%\email{frosini@unisi.it} %
%\and %
%Laboratoire d'Informatique, Algorithmique, Fondements et
%Applications (LIAFA), Universit\'e Denis Diderot 2, place Jussieu
%75251 Paris 5 Cedex 05, France
%\email{Maurice.Nivat@liafa.jussieu.fr}}

\author[Siena]{Andrea Frosini}
\ead{frosini@unisi.it}
\address[Siena]{Dipartimento di Scienze Matematiche ed Informatiche
``Roberto Magari"\\
Universit\`a degli Studi di Siena\\
Pian dei Mantellini 44, 53100, Siena, Italy}%
\author[Paris]{Maurice Nivat}
\ead{Maurice.Nivat@liafa.jussieu.fr}
\address[Paris]{Laboratoire d'Informatique, Algorithmique, Fondements et
Applications (LIAFA), Universit\'e Denis Diderot 2,\\ place
Jussieu 75251 Paris 5 Cedex 05, France }

%\maketitle

%%%%%%%%%%%%%%%%%%%%%%%%%%%%%%%%%%%%%%%%%%%%%%%%%%%%%%%%%%

\begin{abstract}

A binary matrix can be scanned by moving a fixed rectangular
window (sub-matrix) across it, rather like examining it closely
under a microscope. With each viewing, a convenient measurement is
the number of $1s$ visible in the window, which might be thought
of as the {\em luminosity} of the window. The {\em rectangular
scan} of the binary matrix is then the collection of these
luminosities presented in matrix form. We show that, at least in
the technical case of a {\em smooth} $m\times n$ binary matrix, it
can be  reconstructed from its rectangular scan in polynomial time
in the parameters $m$ and $n$, where the degree of the polynomial
depends on the size of the window of inspection. For an arbitrary
binary matrix, we then extend this result by determining the
entries in its rectangular scan that preclude the smoothness of
the matrix.

\end{abstract}
\end{frontmatter}

\section{Introduction and Definitions}

The aim of {\em discrete tomography} is the retrieval of
geometrical information about a physical structure, regarded as a
finite set of points in the integer square lattice $\mathbb{Z}
\times \mathbb{Z}$, from measurements, generically known as {\em
projections}, of the number of atoms in the structure that lie on
lines with fixed scopes (see \cite{HK} for a survey). A common
simplification is to represent a finite physical structure as a
binary matrix, where an entry is $1$ or $0$ according as an atom
is present or absent in the structure at the corresponding point
of the lattice. The challenge is then to reconstruct key features
of the structure from a small number of scans of projections
\cite{Rys}, eventually using some a priori information as
convexity \cite{BDNP} \cite{CD}, and periodicity \cite{FNV}.

Our interest here, following \cite{Nivat}, is to probe the
structure, not with lines of fixed scope, but with their natural
two dimensional analogue, rectangles of fixed scope, much as we
might examine a specimen under a microscope or magnifying glass.
For each position of our rectangular probe, we count the number of
visible atoms, or, in the simplified binary matrix version of the
problem, the number of $1$ in the prescribed rectangular window,
which we term its {\em luminosity}. In the matrix version of the
problem, these measurements can themselves be organized in matrix
form, called the {\em rectangular scan} of the original matrix.
Our first objective is then to furnish a strategy to reconstruct
the original matrix from its rectangular scan. In the sequel, we
will address this problem to as {\em Reconstruction}$(A,p,q)$,
where $A$ is the rectangular scan, and $p$ and $q$ are the
dimension of the rectangular windows. As we also note, our
investigation is closely related to results on tiling by
translation in the integer square lattice discussed in
\cite{Nivat}.

To be more precise, let $M$ be an $m\times n$ integer matrix, and,
for fixed $p$ and $q$, with $1\leq p\leq m, 1\leq q \leq n$,
consider a $p\times q$ window $R_{p,q}$ allowing us to view the
intersection of any $p$ consecutive rows and $q$ consecutive
columns of $M$. Then, the number $R_{p,q}(M)[i,j]$ %of $1s$ in $M$
on view when the top left hand corner of $R_{p,q}$ is positioned
over the $(i,j)$-entry, $M[i,j]$, of $M$, is given by summing all
the entries on view:
$$R_{p,q}(M)[i,j]=\sum_{r=0}^{p-1} \sum_{c=0}^{q-1} M[i+r,j+c],
\quad 1\leq i\leq m-p+1,\quad 1\leq j\leq n-q+1.$$ Thus, we obtain
an $(m-p+1)\times (n-q+1)$ matrix $R_{p,q}(M)$
%with non-negative
%integer entries with entries $R_{p,q}(M)[i,j]$
%as illustrated in Fig.~\ref{smooth1}. We
called the {\em $(p,q)$-rectangular scan} of $M$; when $p$ and $q$
are understood, we write $R(M)=R_{p,q}(M)$, and speak more simply
of the {\em rectangular scan}. (This terminology is a slight
departure from that found in \cite{Nivat}.) In the special case
when $R(M)$ has all entries equal, say $k$, we say that the matrix
$M$ is {\em homogeneous} of {\em degree} $k$, simply
$k$-homogeneous.

%\medskip
%
%\begin{figure}[htd]
%\centerline{\hbox{\psfig{figure=smmooth1.eps,width=2.8in,clip=}}}
%\caption{\label{smooth1} {\em \small A matrix $M$ and its
%$(2,3)$-rectangular scan. The sum of the elements of $M$ inside
%the highlighted rectangle determines the highlighted element of
%$R(M)$.}}
%\end{figure}
%
%\medskip

Furthermore,
%Given any $m\times n$ matrix $A$ and integers $p$ and $q$ with
%$1\leq p \leq m$ and $1\leq q \leq n$,
we define an $(m-p) \times (n-q)$ matrix ${\chi}_{p,q}(M)$
%=({\chi}_{p,q}(M)[i,j])$
by setting, for $1\leq i\leq m-p, \ 1\leq j\leq n-q$:
$${\chi}_{p,q}(M)[i,j]=M[i,j]+M[i+p,j+q]-M[i+p,j]-M[i,j+q].$$
%(see Fig.~\ref{figchi}).
%note that, in the case where A is a {\em binary} matrix, these
%entries take only the values -2, -1, 0, 1, or 2 (see
%Fig.~\ref{figchi}).
As usual, when $p$ and $q$ can be understood without ambiguity, we
suppress them as subscripts. In the event that the matrix
${\chi}(M)$ is a zero matrix, the matrix $M$ is said to be {\em
smooth}. Notice that the homogeneous matrices are {\em properly}
included in the smooth matrices, as shown by the matrix $M$ of
Fig.~\ref{figrev}, which is smooth, and non homogeneous.

\medskip

\begin{figure}[htd]
\centerline{\hbox{\psfig{figure=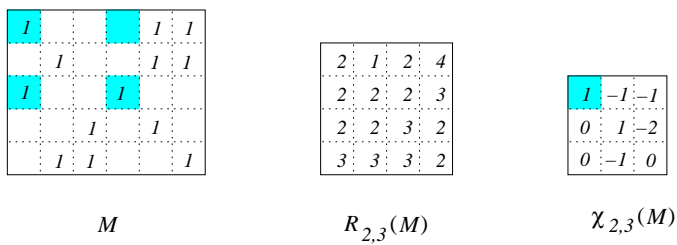,width=3.5in,clip=}}}
\caption{\label{figchi} {\em \small The matrix $M$ and its
corresponding matrices $R_{2,3}(M)$ and $\chi_{2,3}(M)$.}}
\end{figure}

\medskip

Simplifying the rules of the game, thought the paper we will
consider the matrix $M$ (representing a physical structure) as a
{\em binary one}; under this assumption, the rectangular scan
$R(M)$ turns out to be a positive matrix whose values are in the
set $\{0,\dots, p\cdot q \}$, and the matrix $\chi(M)$ turns out
to have values in the set $\{-2,\dots,2\}$ (see
Fig.~\ref{figchi}).

We conclude this introductory section with three observations
which are direct consequences of the given definitions and to
which we shall have frequent recourse in what follows. Since their
proofs are a matter of simple computations, they are omitted.

\begin{lem}\label{rem1}
If $M_1$ and $M_2$ are two $m\times n$ binary matrices, then
$$
R(M_1+M_2)= R(M_1)+R(M_2) \:\: \mbox{ and }\:\:
\chi(M_1+M_2)=\chi(M_1)+\chi(M_2).
$$
\end{lem}

\begin{lem}\label{lem0}
If $M$ is a binary matrix, then
$${\chi}_{1,1}(R_{p,q}(M)) = {\chi}_{p,q}(M). $$
\end{lem}

Thus the rectangular scan $R(M)$ of a binary matrix $M$ already
contains sufficient information to compute $\chi(M)$ and so to
decide whether $M$ is smooth. Notice that, with a certain
terminological inexactitude, we can also say, in the case where
$M$ is smooth, that $R(M)$ is smooth (more precisely, $R(M)$ is
$(1,1)$-smooth, while $M$ itself is $(p,q)$-smooth, as our more
careful statement of the Lemma~\ref{lem0} makes clear).

An appeal to symmetry and induction yields the following
generalization of \cite[Lemma 2.2]{Nivat}.

\begin{lem}\label{lem1}
If $M$ is a smooth matrix then, for any integers $\alpha$ and
$\beta$ such that $1\leq i+\alpha p\leq m$ and $1\leq j+\beta
q\leq n$,
$$
M[i,j]+M[i+\alpha p,j+\beta q]=M[i+\alpha p,j]+M[i,j+\beta q].
$$
\end{lem}

Finally, we say that an entry $M[i,j]$ of the matrix $M$ is {\em
$(p,q)$-invariant} if, for any integer $\alpha$ such that $1\leq
i+\alpha p\leq m$ and $1\leq j+\alpha q\leq n$,
$$ M[i+\alpha p,j+\alpha q]=M[i,j]. $$
If all the entries of $M$ are $(p,q)$-invariant, then $M$ is said
to be $(p,q)$-invariant.

\section{A Decomposition Theorem for Binary Smooth Matrices}

In this section we extend the studies about homogeneous matrices
started in \cite{Nivat} to the class of smooth matrices: first we
furnish a series of simple results which link smoothness and
invariance, then we proceed along a path leading through a
decomposition theorem for smooth matrices to their reconstruction.

\begin{lem}\label{lem2}
If $M$ is a smooth matrix, then each of its elements is
$(p,0)$-invariant or $(0,q)$-invariant.
\end{lem}

\pf Since $M$ is smooth, for each $1\leq i \leq m-p$ and $1\leq j
\leq n-q$, it holds
$$
M[i,j]+M[i+p,j+q]=M[i+p,j]+M[i,j+q].
$$
Let us consider the following three possibilities for the element
$M[i,j]$:

\begin{description}%
\item{$i)$} $M[i,j] \not= M[i+p,j]$: by Lemma~\ref{lem1}, for
$\alpha =1$ and for all $\beta \in \mathbb{Z}$ such that $1\leq
j+\beta q \leq n$, it holds $M[i,j+\beta q]=M[i,j]$ and
$M[i+p,j+\beta q]=M[i+p,j]$, so $M[i,j]$ is $(0,q)$-invariant.

\smallskip

\item{$ii)$} $M[i,j]\not= M[i,j+q]$: by reasoning similarly to
$i)$, we obtain that $M[i,j]$ is $(p,0)$-invariant.

\smallskip

\item{$iii)$} $M[i,j]=M[i,j+q]=M[i+p,j]$: if there exists
$\alpha_0 \in \mathbb{Z}$ such that $M[i+\alpha_0 p, j]\not =
M[i,j]$, again reasoning as in $i)$, we obtain that $M[i,j]$ is
$(0,q)$-invariant.

On the other hand, if for all $1\leq i+\alpha p \leq m$ it holds
that $M[i+\alpha p, j] = M[i,j]$, then $M[i,j]$ is
$(p,0)$-invariant.
\end{description}
\noindent Finally, if $m-p+1 \leq i \leq m$ and $n-q+1 \leq j \leq
n$, a similar reasoning leads again to the thesis. \qed

The reader can check that each entry of the smooth matrix $M$ in
Fig.~\ref{figrev} is $(2,0)$-invariant (the highlighted ones) or
$(0,3)$-invariant. A first decomposition result follows:

\begin{thm}\label{decompM}
A matrix $M$ is smooth if and only if it can be obtained by
summing up a $(p,0)$-invariant matrix $M_1$ and a
$(0,q)$-invariant matrix $M_2$ such that they do not have two
entries $1$ in the same position.
\end{thm}

\pf ($\Rightarrow$) Let $M_1$ and $M_2$ contain the
$(p,0)$-invariant and the $(0,q)$-invariant elements of $M$,
respectively. By Lemma~\ref{lem2}, the thesis is achieved.

($\Leftarrow$) Since $M_1$ is $(p,0)$-invariant, then for each
$1\leq i \leq m-p$, $1\leq j \leq n-q$ it holds
$$
\begin{array}{l}
\chi(M_1)[i,j]=M_1[i,j]+M_1[i+p,j+q]-M_1[i+p,j]-M_1[i,j+q]=\\
\\
=M_1[i,j]+M_1[i,j+q]-M_1[i,j]-M_1[i,j+q]=0
\end{array}
$$

\noindent so, by definition, $M_1$ is smooth. The same result
holds for $M_2$ and, by Lemma~\ref{rem1}, for $M=M_1+M_2$. \qed

\smallskip

\noindent We can go further on by reformulating this last theorem
in terms of the rectangular scans of the matrices $M_1$ and $M_2$:

\begin{lem}\label{lem3}
The following statements hold:

\smallskip

$i)$ if $M$ is $(0,q)$-invariant, then $R(M)$ has
constant rows;%

\smallskip

$ii)$ if $M$ is $(p,0)$-invariant, then $R(M)$ has
constant columns.%
\end{lem}

\pf $i)$ For each $1\leq i \leq m-p+1$ and $1\leq j \leq n-q$, we
prove that $R(M)[i,j]=R(M)[i,j+1]$:
$$
\begin{array}{c}
R(M)[i,j+1]= \sum_{r=0}^{p-1} \sum_{c=0}^{q-1} M[i+r,j+1+c]=\\
\\
=\sum_{r=0}^{p-1} \sum_{c=1}^{q-1} M[i+r,j+c]+\sum_{r=0}^{p-1}
M[i+r,j+q]=
\end{array}
$$
since $M$ is $(0,q)$-invariant
$$
=\sum_{r=0}^{p-1} \sum_{c=1}^{q-1} M[i+r,j+c]+\sum_{r=0}^{p-1}
M[i+r,j]=R(M)[i,j].
$$

$ii)$ The proof is similar to $i)$. \qed

After observing that each matrix having constant rows or columns
is smooth, a direct consequence of Theorem~\ref{decompM} and
Lemma~\ref{lem3} is the following:

\begin{thm}\label{decomp2}
A binary matrix $M$ is smooth if and only if $R(M)$ can be
decomposed into two matrices $R_r$ and $R_c$ having constant rows
and columns, respectively.
\end{thm}

\noindent Fig.~\ref{figrev} shows that the converse of the two
statements of Lemma~\ref{lem3} does not hold in general. However,
we can prove the following weaker version:

\begin{figure}[htd]
\centerline{\hbox{\psfig{figure=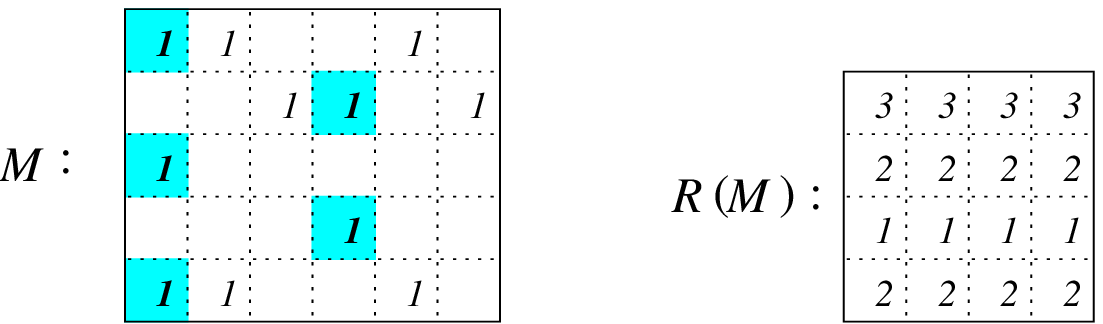,width=3.3in,clip=}}}
\caption{\label{figrev} \em \small A non invariant matrix $M$
whose $(2,3)$-rectangular scan has constant rows.}
\end{figure}

\begin{prop}\label{lem4}
Let $M$ be a binary matrix. The following statements hold:

\smallskip

$i)$ if $R(M)$ has constant columns, then there exists a
$(p,0)$-invariant matrix $M_1$ such that $R(M)=R(M_1)$;

\smallskip

$ii)$ if $R(M)$ has constant rows, then there exists a
$(0,q)$-invariant matrix $M_2$ such that $R(M)=R(M_2)$.
\end{prop}

\pf $i)$ We define the matrix $M_1$ as follows: the first $p$ rows
of $M_1$ are equal to those of $M$, and the other entries of $M_1$
are set according to the desired $(p,0)$-invariance. It is easy to
verify that $R(M_1)=R(M)$.

$ii)$ A definition of $M_2$ similar to that in $i)$ for $M_1$ can
be easily given. \qed

\subsection{Solving $Reconstruction$ $(A,p,q)$ for Smooth Matrices}

A first approach to the general reconstruction problem consists in
the definition of the following algorithm which suites only for a
binary smooth matrix whose $(p,q)$-rectangular scan has constant
rows:

\medskip

\noindent {\sc RecConstRows}$(A,p,q)$

\begin{description}

\item{\bf Input:} an integer matrix $A$ of dimension $m'\times
n'$, having constant rows, and two integers $p$ and $q$.

\item{\bf Output:} a $(0,q)$-invariant matrix $M$, of dimension
$m\times n$, where $m=m'+p-1$ and $n=n'+q-1$, having $A$ as
$(p,q)$-rectangular scan, if it exists, else return FAILURE.

\item{\bf Procedure:}

\item{Step $1$:} create the $m\times n $ matrix $M$ and the vector
$PEnt$ (storing the $P$artial number of $Ent$ries 1 in each row of
$M$) of dimension $m$, to support the computation. Initialize the
entries both of $M$ and $PEnt$ to $0$.

\noindent For each row $1\leq i \leq m'-1$,

\begin{description}

\item{Step $1.1$:} {\bf if} $A[i,1]\leq A[i+1,1]$ {\bf then}

\smallskip

\begin{description}
\item{$\bullet$} $M[i+p,1]=\dots
=M[i+p,A[i+1,1]-A[i,1]+PEnt[i]]=1;$

\smallskip

\item{$\bullet$} $PEnt[i+p]=A[i+1,1]-A[i,1]+PEnt[i].$

\end{description}

\smallskip

{\bf If} $PEnt[i+p]>q$ {\bf then} FAILURE.

\smallskip

\item{Step $1.2$: } {\bf if} $A[i,1]>A[i+1,1]$ and $PEnt[i]\geq
A[i,1]-A[i+1,1]$ {\bf then}

\smallskip

\begin{description}
\item{$\bullet$} $M[i+p,1]=\dots
=M[i+p,A[i+1,1]-A[i,1]+PEnt[i]]=1;$

\smallskip

\item{$\bullet$} $PEnt[i+p]=A[i+1,1]-A[i,1]+PEnt[i].$
\end{description}

\smallskip

\item{Step $1.3$: } {\bf if} $A[i,1]>A[i+1,1]$ and $PEnt[i]<
A[i,1]-A[i+1,1]$ {\bf then}

\smallskip

\begin{description}
\item{$\bullet$} $k=A[i,1]-A[i+1,1]-PEnt[i];$

\smallskip

\item{$\bullet$} for each $i'\leq i$, $i'=(i)mod_p$
\begin{description}

\smallskip

\item{$\centerdot$} $M[i',PEnt[i']+1]=\dots =M[i',PEnt[i']+k]=1;$

\smallskip

\item{$\centerdot$} $PEnt[i']=PEnt[i']+k;$

\smallskip

\item{$\centerdot$} {\bf if} $PEnt[i']>q$ {\bf then} FAILURE.
\end{description}
\end{description}
\end{description}

\smallskip

\item{Step $2$:} let $k=A[1,1]-PEnt[1]- \dots - PEnt[p]$.

For each $1\leq k' \leq k$, search one of the upper leftmost
$p\times q$ positions of $M$, say $(i,j)$, such that, for each
$i'=(i)mod_p$, $1\leq i'\leq m$, it holds $M[i',j]=0$.

{\bf If} such a position does not exist {\bf then} FAILURE,

{\bf else} set all the entries $M[i',j]$ to the value $1$, and
increase $k'$ by one.

\smallskip

\item{Step $3$:} complete the entries of $M$ according to the
$(0,q)$-invariance constraint, and return $M$ as OUTPUT.
\end{description}

As regard the correctness of this reconstruction algorithm, it
relies on the analysis of what stored in $M$ after Step $1$: at
that stage, in fact, the entries in first column of its
rectangular scan $R(M)$ differ from those of $A$ by the same
constant value, without overcoming.

The formal counterpart of what sketched above is in the following
lemmas:

\begin{lem}\label{lem:recrow}
After performing Step $1$ of {\sc RecConstRows}$(A,p,q)$, for each
$1\leq i < m'$, it holds
$$R(M)[i,1]-R(M)[i+1,1]=A[i,1]-A[i+1,1].$$
\end{lem}

\pf Let us first inspect the entries placed in the rows $i$ and
$i+p$ of $M$ during Step $1$ of {\sc RecConstRows}$(A,p,q)$, for a
generic index $1\leq i < m'$:
\begin{description}
\item{if} $A[i,1]\leq A[i+1,1]$, then in row $i+p$ of $M$ are
added $A[i+1,1]-A[i,1]+PEnt[i]$ entries $1$. Since row $i$ of $M$
contains $PEnt[i]$ entries $1$, then, at that step, it holds
\begin{equation}\label{eq1}
R(M)[i,1]-R(M)[i+1]=A[i,1]-A[i+1,i];
\end{equation}

\smallskip

\item{if} $A[i,1]> A[i+1,1]$, and $PEnt[i]\geq A[i,1]-A[i+1,1]$,
then in row $i+p$ of $M$ are added $PEnt[i] - A[i,1]+ A[i+1,1]$
entries $1$, so equation (\ref{eq1}) still holds;

\smallskip

\item{if} $A[i,1]> A[i+1,1]$, and $PEnt[i]< A[i,1]-A[i+1,1]$, then
in row $i$ of $M$ are added $A[i,1]-A[i+1,1]-PEnt[i]$ entries $1$
in addition to the $PEnt[i]$ ones already present, so equation
(\ref{eq1}) is again satisfied.
\end{description}

For each $i<i'<i+p$, Step $1$ eventually changes some entries from
row $i+1$ to row $i+p-1$ of $M$. So, both $R(M)[i,1]$ and
$R(M)[i+1,1]$ increase their value of the same amount, without
compromising the validity of equation (\ref{eq1}).

Finally, if $i'\geq i+p$, then Step $1.3$ may modify the values of
$R(M)[i,1]$ and $R(M)[i+1,1]$, but again of the same amount, since
the (eventually) added entries $1$ respect the $(p,0)$-invariance
in the rows of index less than $i'$, so equation (\ref{eq1})
definitively holds, and we obtain the thesis. \qed

\begin{lem}\label{lem:maxim}
Let us consider the vector $PEnt$ as updated at the end of Step
$1$ of {\sc RecConstRows}$(A,p,q)$. For each matrix $M$ such that
$R(M)=A$, and for each $0\leq i \leq m$, it holds
$$
M[i,1]+\dots +M[i,p]\geq PEnt[i].
$$
\end{lem}

\pf By contradiction, we assume that there exists an index $1\leq
i \leq m$ and a matrix $M'$ such that
%\begin{equation}\label{eq2}
$$
M[i,1]+\dots +M[i,p]+k= PEnt[i],
$$
%\end{equation}

with $R(M)=A$ and $k>0$. By Lemma~\ref{lem:recrow}, the same
equation holds for each row $i'=(i)mod_p$.

Let $i_0$ be the first index such that
\begin{description}
\item{$\bullet$} $i_0=(i)mod_p$;

\smallskip

\item{$\bullet$} for each $i'=(i)mod_p$, it holds $PEnt[i_0]\leq
PEnt[i']$.
\end{description}

If $i_0\leq m-p$, then the minimality of the value of $PEnt[i_0]$
assures that $PEnt[i_0]=0$ before Step $1$ reached row $i_0$, and,
consequently, that $A[i_0,1]\leq A[i_0+1,1]$.

As soon as Step $1.1$ reaches the row index $i_0$, it eventually
increases the value $PEnt[i_0+p]$, leaving unchanged that of
$PEnt[i_0]$. Again the minimality of $PEnt[i_0]$ assures that no
changes will be performed to the value of $PEnt[i_0]$ till the end
of Step $1$.

Hence, the assumption $ M[i_0,1]+\dots +M[i_0,p]+k= PEnt[i_0] $,
with $k>0$, generates a contradiction.

If $i_0>m-p$, then a similar argument holds, and so we get the
thesis. \qed

\medskip

Now, also Step $2$ of {\sc RecConstRows} $(A,p,q)$ can be better
understood: the $k$ elements $0$ which change their value to $1$
and which are added in the upper leftmost $p\times q$ positions of
$M$, fill the gap among $R(M)[i,1]$ and $A[i,1]$, so that the
output matrix $M$ has the desired property $R(M)=A$.

\begin{cor}
Each row $i$ of a matrix $M$ having $A$ as rectangular scan
contains at least $PEnt[i]$ elements which are $(p,0)$-invariant,
and not $(0,q)$-invariant.
\end{cor}

A procedure which reconstructs a smooth matrix whose
$(p,q)$-rectangular scan $A$ has constant columns, say {\sc
RecConstCols}$(A,p,q)$, can be easily inferred from {\sc
RecConstRows}, so, in the sequel, we will consider it as already
defined.

From Lemmas \ref{lem:recrow} and \ref{lem:maxim}, it is
straightforward that

\begin{thm}\label{comprec}
The problem $Reconstruction(A,p,q)$ can be solved in $O(m\: n)$,
when $A$ has constant rows or constant columns.
\end{thm}

\begin{exmp}\label{ex0}
{\small Let us follow the computation {\sc RecConstRows}$(A,3,4)$,
with $A$ depicted in Fig.\ref{figxmp0}.

\begin{figure}[htd]
\centerline{\hbox{\psfig{figure=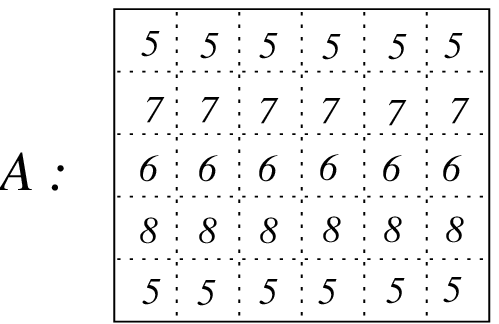,width=1.5in,clip=}}}
\caption{\label{figxmp0} The matrix $A$ of Example~\ref{ex0}.}
\end{figure}

Step $1$: the matrix $M$ is created, and its first four columns,
together with the vector $PEnt$, are modified as shown in
Fig~\ref{figexmp41}. More precisely,
\begin{description}

\item{$A[1,1]+2=A[1,2]$} requires Step $1.1$ to place two entries
$1$ in (the leftmost positions of) row $4$, Fig~\ref{figexmp41},
(a);

\smallskip

\item{$A[1,2]=A[1,3]+1$} requires Step $1.3$ to place one entry
$1$ in row $2$, Fig~\ref{figexmp41}, (b);

\smallskip

\item{$A[1,3]+2=A[1,4]$} requires Step $1.1$ to place two entries
$1$ in row $6$, Fig~\ref{figexmp41}, (c);

\smallskip

\item{$A[1,1]=A[1,2]+3$} requires Step $1.3$ to add one entry $1$
both in row $1$ and in row $4$, Fig~\ref{figexmp41}, (d).

\end{description}

\medskip

\begin{figure}[htd]
\centerline{\hbox{\psfig{figure=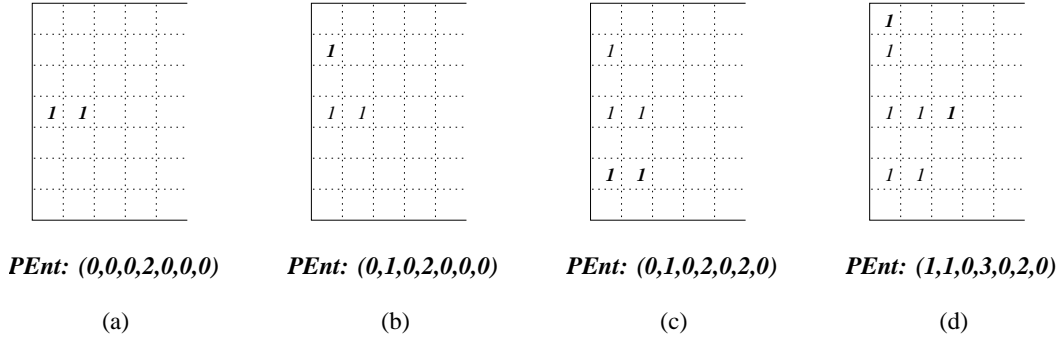,width=5.5in,clip=}}}
\caption{\label{figexmp41} \em \small Step $1$ of {\sc
RecConstRows}$(A,3,4)$.}
\end{figure}

Step $2$ places the remaining $k=A[1,1]-PEnt[1]-PEnt[2]-PEnt[3]=3$
entries $1$ in the upper leftmost $p\times q$ submatrix of $M$,
and propagates them according to the $(3,0)$-invariance, paying
attention that no collisions occur (see Fig.~\ref{figexmp42}
$(a)$).

Step $3$ completes $M$ according to the $(0,4)$-invariance, giving
the final solution depicted in Fig.~\ref{figexmp42} $(b)$.

\medskip

\begin{figure}[htd]
\centerline{\hbox{\psfig{figure=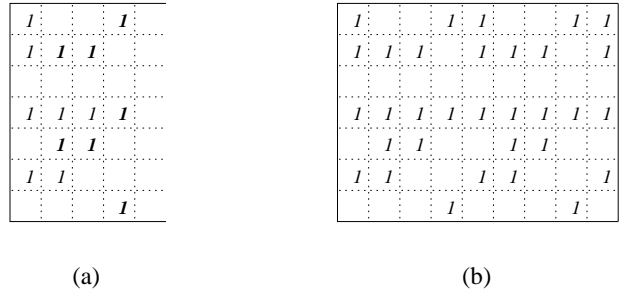,width=3.2in,clip=}}}
\caption{\label{figexmp42} \em \small Step $2$ and Step $3$ of
{\sc RecConstRows}$(A,3,4)$.}
\end{figure}

\medskip
}
\end{exmp}

\noindent {\em The general reconstruction algorithm}

\medskip

Theorem~\ref{decomp2} allows one to foresee the use of the
procedures {\sc RecConstRows} and {\sc RecConstCols} to solve {\em
Reconstruction} $(A,p,q)$, when $A$ is $(1,1)$-smooth:
%, or equivalently, when its solutions are $(p,q)$-smooth.
the algorithm at first will split the matrix $A$ into two parts
having constant rows and columns, respectively, then it will apply
to each of them the appropriate reconstruction procedure, and
finally it will merge the two outputs. Performing the merging
stage a conflict occurs when the same position in two output
matrices has value $1$. To prevent it small refinements to the
outputs of {\sc RecConstRows} and {\sc RecConstCols} will be
required.

So, let us start by showing in the next lemma a quick way of
finding all the possible decompositions of a $(1,1)$-smooth matrix
into two parts having constant rows and columns, respectively.

\begin{lem}\label{lem5}
Let $A$ be a $m\times n$ integer matrix. If $A$ is $(1,1)$-smooth,
then it admits $k+1$ different decompositions into two matrices
having constant rows and columns, with $k$ being the minimum among
all the elements of $A$.
\end{lem}

\pf The thesis is achieved by defining a procedure which gives as
output a complete list of couples of matrices $(A_r^t,A_c^t)$,
with $0\leq t \leq k$, each of them representing a decomposition
of $A$ into two parts having constant rows and columns, and
successively, by proving its correctness:

\medskip

\noindent {\sc Decompose} $(A)$

\begin{description}
\item{\bf Input:} an integer $m\times n$ matrix $A$.

\item{\bf Output:} a sequence of different couples of matrices
$(A_r^0,A_c^0), \dots , (A_r^{k},A_c^{k})$, with $k$ being the
minimum element of $A$, such that, for each $1\leq t \leq k$,
$A_r^t$ has constant rows, $A_c^t$ has constant columns, and
$A_r^t+A_c^t=A$. If such a sequence does not exist, then return
FAILURE.

\item{\bf Procedure:}

\item{Step 1:} initialize all the elements of two $m\times n$
matrices $A_c$ and $A_r$ to the value $0$. Let $k$ be the minimum
among the entries of $A$. From each element of $A$, subtract the
value $k$ and store the result in $A_c$;

\item{Step 2:} for each $1\leq i \leq m$

\begin{description}

\item{Step 2.1:} compute
$$
k_i=min_{j}\{A_c[i,j]\: : \: 1\leq j \leq n\};
$$

\item{Step 2.2:} subtract the value $k_i$ from each element of
$A_c$;

\item{Step 2.3:} set all the elements of row $i$ of $A_r$ to the
value $k_i$;

\end{description}

\item{Step 3:} {\bf if} the matrix $A_c$ has not constant columns
{\bf then} FAILURE

{\bf else} for each $0\leq t \leq k$, create the matrices $A_r^t$
and $A_c^t$ such that
$$
A_r^t[i,j]=A_r[i,j]+t \hspace{20pt}\mbox{ and } \hspace{20pt}
A_c^t[i,j]=A_c[i,j]+k-t,
$$
with $1\leq i \leq m$ and $1\leq j \leq n$.

Give the sequence $(A_r^0,A_c^0),\dots ,(A_r^k,A_c^k)$ as OUTPUT.
\end{description}

Example~\ref{ex1} shows a run of the algorithm. By construction,
each couple $(A_r^t,A_c^t)$ is a decomposition of $A$, and
furthermore, $A_r^t$ has constant rows.

What remains to prove is that the matrix $A_c$ updated at the end
of Step $2$ has constant columns (and, consequently, the same hold
for all the matrices $A_c^t$). Let us denote by $r_i$ the common
value of the elements of the $i$-th row of $A_r$, and let us
proceed by contradiction, assuming that $A_c$ has not constant
columns. Since $A$ is the sum of a column constant and a row
constant matrix, and for all $1\leq i \leq m-p+1$ and $1\leq j
\leq n-q+1$, it holds $A[i,j]=r_i + A_c[i,j]+k$, then $A_c$ is
also the sum of a column constant matrix, and a row constant
matrix, this last having at least one row, say $i_0$, whose
elements have value $k_{i_0}'\not =0$.

This situation generates an absurd, since $k_{i_0}$ computed in
Step $2.1$ turns out no longer to be the minimum of row $i_0$ in
$A_c$, updated to that step.

Since a matrix having constant rows (resp. columns) cannot be
obtained as sum of a matrix having constant rows and a matrix
having constant columns unless the latter is a constant matrix,
then the $k+1$ decompositions listed by the algorithm are all the
possible ones. \qed

\begin{exmp}\label{ex1}
{\small Let us follow the steps of the procedure {\sc
Decompose}$(A)$, with the matrix $A$ depicted in
Fig.\ref{smooth21}.

\begin{figure}[htd]
\centerline{\hbox{\psfig{figure=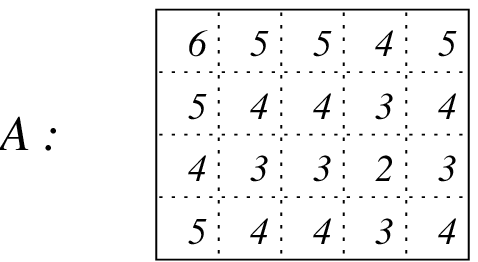,width=1.5in}}}
\caption{\label{smooth21} \em \small The $(1,1)$-smooth matrix $A$
of Example~\ref{ex1}.}
\end{figure}

Step $1$: we subtract from all the elements of $A$, the value
$k=2$, i.e. its minimum element, and we store the obtained result
in the matrix $A_c$.

Step $2$: for each $1\leq i \leq m-p+1$, we find the minimum value
$k_i$ among the elements of row $i$ of $A_c$ (Step $2.1$), we
subtract it from all these elements (Step $2.2$), and finally, we
set the elements in row $i$ of $A_r$ to the value $k_i$ (Step
$2.3$). In our case, the minimums are $k_1=2$, $k_2=1$, $k_3=0$
and $k_4=1$.

\begin{figure}[htd]
\centerline{\hbox{\psfig{figure=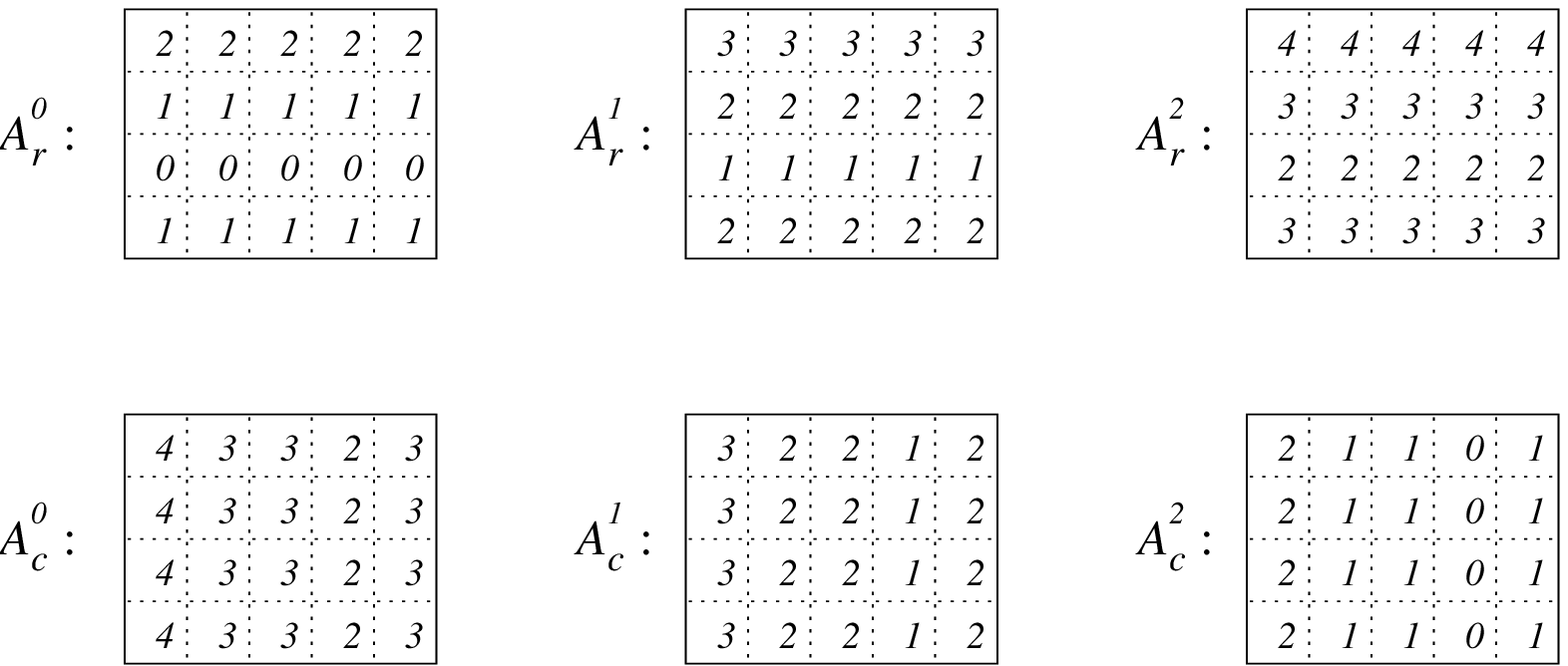,width=4.8in,clip=}}}
\caption{\label{smooth22} \em \small The three decompositions of
the matrix $A$.}
\end{figure}

Step $3$: the matrix $A_c$ updated at the end of Step $2$ has
constant columns, so the three different decompositions of $A$ can
be computed and listed.

\noindent The output is depicted in Fig.~\ref{smooth22}.}
\end{exmp}

Now we are finally able to define a general procedure which solves
the problem {\em Reconstruction} $(A,p,q)$, when $A$ is
$(1,1)$-smooth:

\medskip

\noindent {\sc RecSmooth} $(A,p,q)$

\begin{description}

\item{\bf Input:} an integer $(1,1)$-smooth matrix $A$ of
dimension $m'\times n'$ and two integers $p$ and $q$.

\smallskip

\item{\bf Output:} a binary matrix $M$ of dimension $m\times n$,
with $m=m'+p-1$ and $n=n'+q-1$, having $A$ as $(p,q)$-rectangular
scan, if it exists, else return FAILURE.

\smallskip

\item{\bf Procedure:}

\item{Step $1$:} run {\sc Decompose}$(A)$, and let $(A^0_r,A^0_c),
\dots ,(A^k_r,A^k_c)$ be its output.

Set $t=0$;

\smallskip

\begin{description}
\item{Step $1.1$:} run Step $1$ of {\sc RecConstRows}
$(A^t_r,p,q)$.

Let $PEnt_{row}=PEnt$ and $k_{row}=k$, with $PEnt$ and $k$ updated
at the end of the step. Define $PEnt_{\equiv p}$ to be the vector
having $p$ elements, and such that:
$$
PEnt_{\equiv p}[i]=max\{PEnt_{row}[i']: i'=(i)mod_p\},
$$
with $1\leq i \leq p$, and $1\leq i' \leq m$.

\smallskip

\item{Step $1.2$:} run Step $1$ of {\sc
RecConstCols}$(A^t_c,p,q)$.

Let $PEnt_{col}=PEnt$ and $k_{col}=k$, with $PEnt$ and $k$ updated
at the end of the step, and let $PEnt_{\equiv q}$ be the vector
having $q$ entries, and such that:
$$
PEnt_{\equiv q}[j]=max\{PEnt_{col}[j']: j'=(j)mod_q\},
$$
with $1\leq j \leq q$, and $1\leq j' \leq n$;

\smallskip

\item{Step $1.3$:} among all the possible $p\times q$ matrices
whose entries are in $\{P,Q,1,0\}$, choose one, say $W$, such
that:

$i)$ the number of entries $Q$ in its $i$-th row is $PEnt_{\equiv
p}[i]$.

\smallskip

$ii)$ the number of entries $P$ in its $j$-th column is
$PEnt_{\equiv q}[j]$;

\smallskip

$iii)$ the number of entries $1$ is $k_{row}+k_{col}$.

\smallskip

{\bf If} $W$ does not exist and $t\not=k$ {\bf then} set $t=t+1$,
and return to Step~$1.1$.

{\bf If} $W$ does not exist and $t=k$ {\bf then} FAILURE;
\end{description}

\medskip

\item{Step $2$:} create the $m\times n$ matrix $M$, and initialize
its entries as follows:

\smallskip

\begin{description}
\item{Step $2.1$:} for each $0\leq i \leq m$ and for each $0\leq j
\leq q$,

{\bf if} $PEnt_{row}[i]\not=0$ and $W[i',j]=Q$, with $i'=(i)mod_p$
{\bf then} set both $M[i,j]=Q$ and $PEnt_{row}[i]=
PEnt_{row}[i]-1$;

\smallskip

\item{Step $2.2$:} for each $0\leq j \leq n$ and for each $0\leq i
\leq p$,

{\bf if} $PEnt_{col}[j]\not=0$ and $W[i,j']=P$, with $j'=(j)mod_q$
{\bf then} set both $M[i,j]=P$ and
$PEnt_{col}[j]=PEnt_{col}[j]-1$;

\smallskip

\item{Step $2.3$:} for each $0\leq i \leq p$ and for each $0\leq j
\leq q$,

{\bf if} $W[i,j]=1$, {\bf then} set $M[i,j]=1$;

\smallskip

\item{Step $2.4$:} fill the matrix $M$ imposing the
$(0,q)$-invariance of its entries $Q$ and $1$, and the
$(p,0)$-invariance of its entries $P$ and $1$;
\end{description}

\medskip

\item{Step $3$:} change the values $P$ and $Q$ to $1$, and set the
remaining entries of $M$ with the value $0$; finally, give $M$ as
output.

\end{description}

\medskip

\begin{thm}
The problem $Reconstruction$ $(A,p,q)$, with $A$ being
$(1,1)$-smooth, admits a solution if and only if {\sc RecSmooth}
$(A,p,q)$ does not return FAILURE.
\end{thm}
\pf $(\Rightarrow)$ Let $M$ be a solution of $Reconstruction$
$(A,p,q)$, and let us assume that $M=M_1+M_2$, with $M_1$ and
$M_2$ being $(0,q)$-invariant and $(p,0)$-invariant, respectively.

Let $(A^t_r,A^t_c)$ be one of the decomposition of $A$
% obtained by {\sc Decompose} $(A)$,
such that $R(M_1)=A^t_r$ and $R(M_2)=A^t_c$.

Lemma~\ref{lem:maxim} implies that, for each $1 \leq i \leq m$,
the value $PEnt_{row}[i]$ indicates the minimum number of elements
of $M$ which are $(0,q)$-invariant and not $(p,0)$-invariant, and
which lie in the first $q$ columns of the solution; a symmetrical
property holds for $PEnt_{col}$.

Let us construct a $p\times q$ matrix $W'$ as follows:

\begin{description}

\item{-} for each $1\leq i \leq m$ and $1\leq j \leq q$, if
$M[i,j]=1$ is $(0,q)$-invariant and not $(p,0)$-invariant, then
set $W'[i',j]=Q$, with $i'=(i)mod_p$;

\smallskip

\item{-} for each $1\leq i \leq p$ and $1\leq j \leq n$, if
$M[i,j]=1$ is $(p,0)$-invariant and not $(0,q)$-invariant, then
set $W'[i',j]=P$, with $j'=(j)mod_q$;

\smallskip

\item{-} for each $1\leq i \leq p$ and $1\leq j \leq q$, if
$M[i,j]=1$ is $(p,0)$-invariant and $(0,q)$-invariant, then set
$W'[i,j]=1$.

\end{description}

Obviously, by definition of invariance, in $W'$ there are no
positions which are first set to a value and then modified to
another. So, the existence of matrix $W'$ implies that of a matrix
$W$ satisfying the constraints imposed in Step $1.3$.

%No further possibilities of failure are left to {\sc RecSmooth}
%$(A,p,q)$.

$(\Leftarrow)$ Immediate. \qed

\medskip

\begin{exmp}\label{ex3}
Let us describe a run of {\sc RecSmooth} $(A,3,3)$ starting from
the decomposition of $A$ into the couple of matrices depicted in
Fig.~\ref{sm8}.

\begin{figure}[htd]
\centerline{\hbox{\psfig{figure=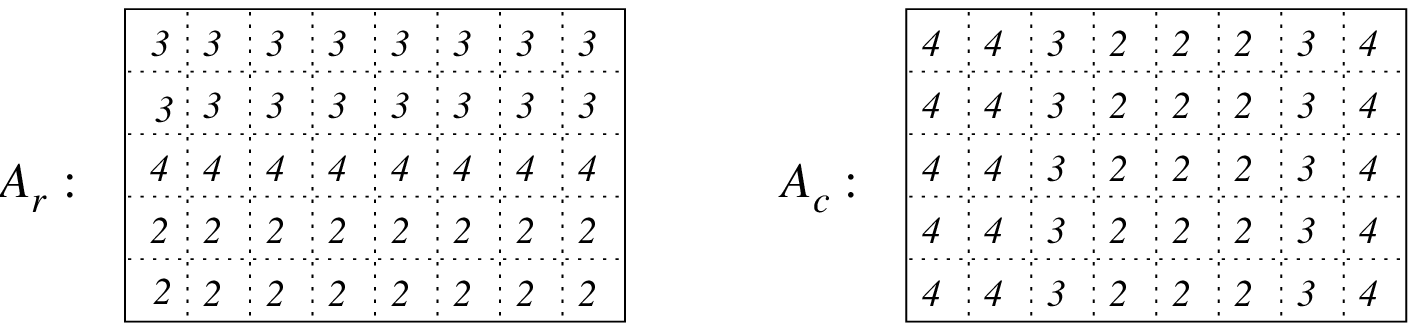,width=4.5in,clip=}}}
\caption{\label{sm8} \em \small The decomposition of the
rectangular scan $A$ used in Example~\ref{ex3}.}
\end{figure}

Step $1.1$ produces the vectors
$$
PEnt_{row}=(0,0,2,0,1,0,0) \;\; \mbox{ and } \;\; PEnt_{\equiv
p}=(0,1,2),
$$
while Step $1.2$ produces the vectors
$$
PEnt_{col}=(0,1,1,0,0,0,0,0,1,1) \;\; \mbox{ and } \;\;
PEnt_{\equiv q}=(1,1,1).
$$

%We reconstruct a matrix $M=M_1+M_2$ such that $R(M)=A$,
%$R(M_1)=A_r$ and $R(M_2)=A_c$.
%
%We start to compute and list all the solutions of {\sc
%RecConstRows}$(A_r,3,3)$ and {\sc RecConstCols}$(A_c,3,3)$, and we
%arrange them in two sequences $s_1$ and $s_2$ as requested in Step
%$1.2$. Then, following Step $1.3$, the elements of $s_1$ and $s_2$
%are summed up in all possible ways until the binary matrix $M$ is
%obtained.
%
%\medskip
%
%\begin{figure}[htd]
%\centerline{\hbox{\psfig{figure=smooth9.eps,width=5in,clip=}}}
%\caption{\label{figexmp21} \em \small Two non-disjoint solutions
%$M_1$ and $M_2$ which are $(0,3)$-invariant and $(3,0)$-invariant,
%as stated in Lemma~\ref{lem4}.}
%\end{figure}
%
%\medskip

\begin{figure}[htd]
\centerline{\hbox{\psfig{figure=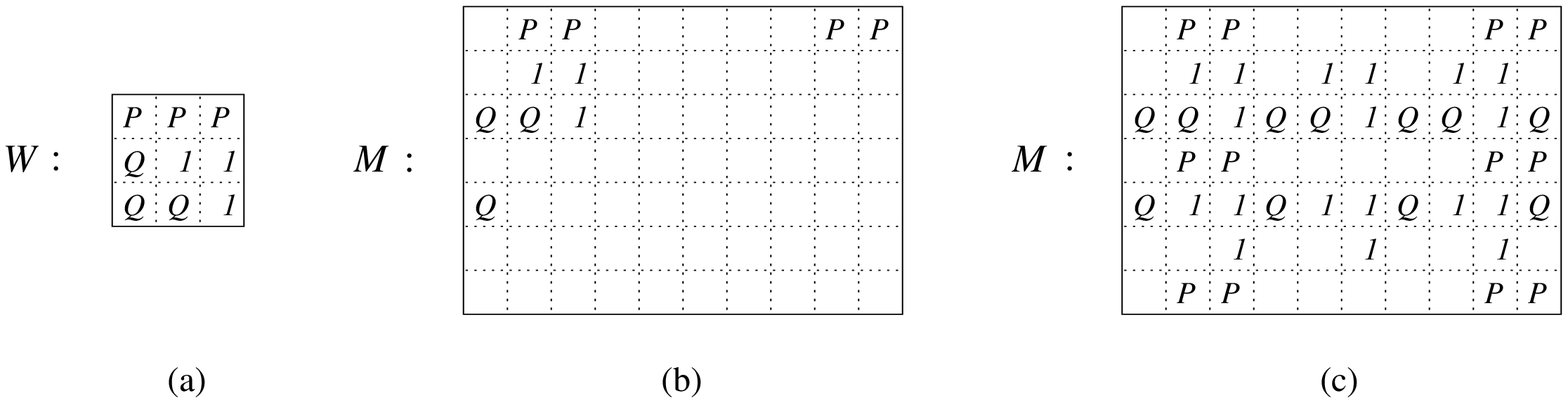,width=6in,clip=}}}
\caption{\label{newfigexmp2} \em \small The matrices created in
successive steps of {\sc RecSmooth} $(A,3,3)$.}
\end{figure}

Among the matrices which are compatible with the requirements of
Step $1.3$ we choose that depicted in Fig.~\ref{newfigexmp2}, (a).

Steps $2.1$, $2.2$ and $2.3$ produce the matrix $M$ in
Fig.~\ref{newfigexmp2}, (b) (notice that the entries $P$ are
$(p,0)$-invariant, while the entries $Q$ are $(0,q)$-invariant as
desired). Finally, Step $2.4$ produces the matrix in
Fig.~\ref{newfigexmp2}, (c), and, consequently, Step $3$ the
output.

\end{exmp}

The following theorem holds:

\begin{thm}\label{teo1}
The computational complexity of {\sc RecSmooth}$(A,p,q)$ is
polynomial in $m$ and $n$.
\end{thm}

\pf We obtain the thesis by analyzing the complexity of each step
of {\sc RecSmooth}:

\begin{description}

\item{Step $1$:} the procedure {\sc Decompose} clearly acts in
$O(m \: n)$ time (remind that $k\leq p \cdot q$ is the minimum
among the elements of $A$);

\item{Step $1.1$} and Step $1.2$ are repeated at most $k$ times,
and, each time, they ask for a run of {\sc RecConstRows} and of
{\sc RecConstCols} which are both performed in $O(m\: n)$. The
computation of $PEnt_{\equiv p}$ and $PEnt_{\equiv q}$ does not
increase the complexity of these two steps.

\item{Step $1.3$} is carried on in constant time with respect to
$m$ and $n$.

\item{Steps $2$} and $3$ require $O(m\: n)$ to create matrix $M$.
\end{description}

Hence, the total amount of time is $O(m\: n)$. \qed

\begin{rem}
We are aware that Step $1.3$ of {\sc RecSmooth} $(A,p,q)$, i.e.
the search of the matrix $W$, can be carried on in a smarter way,
but this will bring no effective contribution to the decreasing of
the computational complexity of the reconstruction, and, on the
other hand, it will add new lemmas and proofs to the current
section.
\end{rem}

The first part of the paper devoted to the analysis and the
reconstruction of smooth matrices is now completed.

\section{Solving $Reconstruction$ $(A,p,q)$: final challenge}

This last section concerns the matrices which are not smooth: in
particular, for each non smooth matrix $M$, we consider the matrix
$\chi(M)$ and we define a polynomial time algorithm which lists
all the matrices consistent with it. Finally we will integrate it
with the algorithm for reconstructing a smooth matrix defined in
the previous section, and we will achieve the solution of the
general reconstruction problem.

Unfortunately, the definitions introduced up to now are not
specific enough to describe these further studies, and a final
effort is required to the reader: what follows has the appearance
of a stand alone part inside this section, but the feeling of a
final possible usage will never be frustrated.

\begin{figure}[htd]
\centerline{\hbox{\psfig{figure=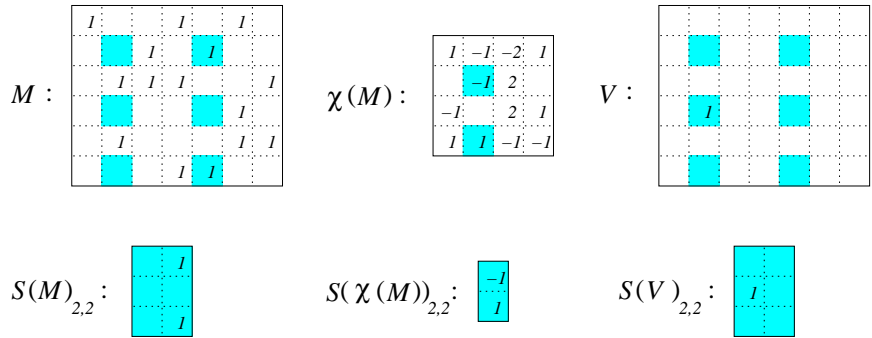,width=4.5in,clip=}}}
\caption{\label{figgrid} \em \small The subgrids of the matrices
$M$ and $\chi(M)$ with respect to the position $(2,2)$. Matrix $V$
is one of the possible valuations of $S(\chi(M))_{2,2}$.}
\end{figure}

Hence, let $a$ and $b$ be two indexes such that $1\leq a \leq p$,
$1\leq b \leq q$, and $A$ be an integer $m\times n$ matrix. We
define the {\it $(a,b)$-subgrid} of $A$ to be the submatrix
$$
S(A)_{a,b}[i,j]=A[a+(i-1)\: p,b+(j-1)\: q]
$$
with $1 \leq a+(i-1)p \leq m$ and $1\leq b+(j-1)q \leq n$
%where $\frac{1-a+p}{p} \leq i \leq \frac{m-a+p}{p}$, and
%$\frac{1-b+q}{q} \leq j \leq \frac{m-b+q}{q}$
(see Fig.~\ref{figgrid}).
%$a+(i-1)\: p \leq m$, and $b+(j-1)\: q\leq n$
%(so $S(A)_{a,b}$ turns out to have dimension $\lceil (m-a+1)/p
%\rceil \times \lceil (n-b+1)/q \rceil$).

If we consider again a binary matrix $M$, by definition it holds
that
$$
\begin{array}{l}
\chi(M)[a+(i-1)\: p,b+(j-1)\: q]=S(\chi(M))_{a,b}[i,j]=\\
\\
=S(M)_{a,b}[i,j]+S(M)_{a,b}[i+1,j+1]-S(M)_{a,b}[i+1,j]-
S(M)_{a,b}[i,j+1].
\\
\end{array}
$$
The binary matrix $V$ of dimension ${m\times n}$ is said to be a
{\it valuation} of $S(\chi(M))_{a,b}$ if, for each $1\leq i \leq
m$, $1\leq j \leq n$,
\begin{description}
\item{- } if $i\not = (a) mod_p$ and $j\not = (b) mod_q$ then
$V[i,j]=0$;

\smallskip

\item{- } $S(\chi(M))_{a,b}=S(\chi(V))_{a,b}$ (see
Fig.~\ref{figgrid}).
\end{description}

The notion of valuation extends to the whole matrix $\chi(M)$ as
the union of the valuations of all its subgrids.

\begin{prop}\label{lem8}
Let $S(\chi (M))_{a,b}$ and $S(\chi (M))_{a',b'}$ be two subgrids
whose valuations are $V$ and $V'$, respectively. If $a\not= a'$ or
$b\not=b'$, then for each $1\leq i \leq m$ and $1\leq j \leq n$,
$V[i,j]=1$ implies $V'[i,j]=0$ .
\end{prop}

%So, the positions having value $1$ of two valuations of different
%subgrids are disjoint.

\begin{lem}\label{lem6}
Let $V$ be a valuation of $S(\chi(M))_{a,b}$, and let $i_0$ be a
row [column] of $S(V)_{a,b}$ having all the elements equal to $1$.
The matrix $V'$ such that $S(V')_{a,b}$ is equal to $S(V)_{a,b}$
except in the elements of the row [column] $i_0$ which are all set
to $0$, is again a valuation of $S(\chi(M))_{a,b}$.
\end{lem}

%A similar result obviously holds if we consider a column of
%$S(V)_{a,b}$.

If $V$ and $V'$ are two valuations as in Lemma~\ref{lem6}, then we
say that the valuation $V$ is {\it greater than} the valuation
$V'$. This relation can be easily extended to a finite partial
order on the valuations of the subgrids of $\chi(M)$.

\begin{lem}\label{rem_uniq}
Let $1\leq i \leq m-p$ and $1\leq j \leq n-q$. If
$\chi(M)[i,j]=2$, then $M[i,j]=M[i+p,j+q]=1$, and $M[i+p,j]=
M[i,j+q]=0$.
\end{lem}

The proof is immediate. A symmetric result holds if $\chi(M)[i,j]$
has value $-2$.

The following lemma turns out to be crucial in this section. A
deeper analysis of what it states could furnish material for
further studies:

\begin{lem}\label{lem7}
Given a binary matrix $M$, for each couple of integers $1\leq a
\leq p$, $1\leq b\leq q$, the number of minimal elements in the
partial ordering of the valuations of $S(\chi (M))_{a,b}$ is
polynomial with respect to the dimensions $m$ and $n$ of $M$.
Furthermore, each minimal element can be reconstructed in
polynomial time with respect to them.
\end{lem}

\pf

%%%%%%%%%%%%%%%%%%%%%%%%%%  INIZIO PROVA %%%%%%%%%%%%%%%%%

Let $S(\chi(M))_{a,b}$ have dimension $m'\times n'$. We order the
(positions of the) non zero elements of $S(\chi(M))_{a,b}$
according to the numbering of its columns (from left to right),
and, in the same column, according to the numbering of its rows
(from up to bottom), and let $p_1, \dots ,p_t$ be the obtained
sequence. We prove the thesis by induction on the number $t$ of
elements of the sequence, i.e. we prove that the addition of new
nonzero elements in $S(\chi(M))_{a,b}$ does not increase ``too
much" the number of its possible minimal valuations.

We first observe that, by Remark~\ref{rem_uniq}, the presence of
entries $2$ or $-2$ in $S(\chi(M))_{a,b}$ does not increase the
number of minimal valuations, so we are allowed to focus our
attention exactly on the elements of value $1$ or $-1$. As one can
expect, the symmetry of the two cases allows us to show the
details of one (in particular when the element to add has value
$1$), and let the reader infer the other:

\begin{description}

\item{\bf Base $t=1$:} if $p_1=(i,j)$ and $S(\chi(M))_{a,b}
[i,j]=1$, then the four possible valuations of $S(\chi(M))_{a,b}$
are depicted in Fig.~\ref{figP1}. Among them, only $S(V_1)_{a,b}$
and $S(V_2)_{a,b}$ are minimal: they can be reached both from
$S(V_3)_{a,b}$ and from $S(V_4)_{a,b}$ by deleting the rows or the
columns entirely filled with entries $1$, as stated in
Lemma~\ref{lem6}.

\begin{figure}
\centerline{\hbox{\psfig{figure=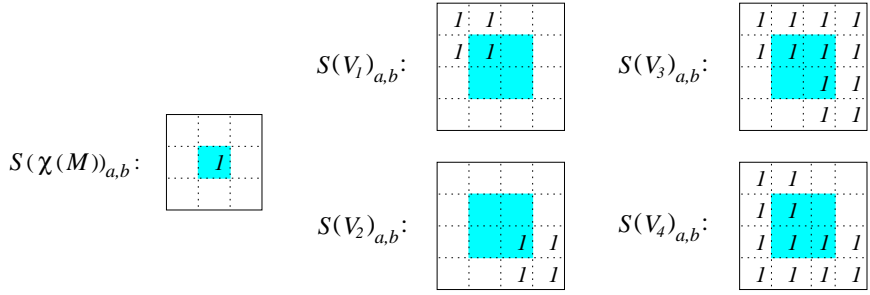,width=4.5in,clip=}}}
\caption{\label{figP1} \em \small
%The four valuations when $t=1$,
%$p_1=(2,2)$, and $S(\chi (M))_{a,b}[2,2]=1$.
The four valuations of the point $p_1=(2,2)$ of value $1$.}
\end{figure}

%%%%%%%%%%%%   SI CONSIDERA UN SOLO PUNTO DA AGGIUNGERE

\item{\bf Step $t\rightarrow t+1$}: let $S(\chi (M))_{a,b}$ have a
sequence $p_1,\dots , p_{t+1}$ of nonzero points, with
$p_{t+1}=(i,j)$, and $S(\chi(M))_{a,b}[i,j]=1$. Let $V$ be a
valuation of the first $p_1,\dots p_t$ points in
$S(\chi(M))_{a,b}$. It is straightforward that, for all $1 \leq i'
\leq m'+1$ and $j< j' \leq n'+1$ it holds
$S(V)_{a,b}[i',j+1]=S(V)_{a,b}[i',j']$.

Hereafter, we show {\em all} the possible ways of extending $V$ to
the valuation $V'$ which includes the point $p_{t+1}$. Some
pictures are supplied in order to make the different cases
transparent.

Let us call
\begin{description}

\item{\em $0$-row}: a row of $S(V)_{a,b}$ whose elements have all
value $0$;

\item{\em $*1$-row}: a row of $S(V)_{a,b}$ whose element in column
$j+1$ has value $1$;

\item{\em $*0$-row}: a row of $S(V)_{a,b}$ which is neither
$0$-row nor $1$-row.
\end{description}

We examine all the possible configurations of $S(V)_{a,b}$, and
for each of them we indicate the desired extension to the {\em
minimal} valuation $S(V')_{a,b}$:

\begin{description}

\item{$i)$} all the rows from $1$ to $i$ are $0$-rows or $*1$-rows
(see Fig.~\ref{sm40}, $(a)$).

We define the valuation $S(V')_{a,b}$ as follows: for each $1\leq
i' \leq i$,
\begin{description}
\item{if} row $i'$ is a $0$-row, then change from $0$ to $1$ the
value of each entry of $S(V)_{a,b}$ in position $(i',j')$, with
$1\leq j' \leq j$, so that it becomes a $*0$-row;

\item{if} row $i'$ is a $*1$-row, then change from $1$ to $0$ the
value of each entry of $S(V)_{a,b}$ in position $(i',j')$, with
$j+1 \leq j' \leq n'+q-1$, so that it becomes a $*0$-row or a
$0$-row. If a $0$-row is created, then discard the obtained
valuation $S(V')_{a,b}$, since it has been already obtained in a
previous step (easy check);
\end{description}

\item{$ii)$} there exists a $*0$-row $i'$, with $1\leq i' \leq i$
(see Fig.~\ref{sm40}, $(b)$).

No changes in the first $i$ rows of $S(V)_{a,b}$ allow the
insertion of the new point $p_{t+1}$;

\item{$iii)$} all the rows from $i+1$ to $m'+p-1$ are $0$-rows or
$*0$-rows (see Fig.~\ref{sm40}, $(c)$).

We define the valuation $S(V')_{a,b}$ as follows: for each
$i+1\leq i' \leq m'+p-1$, change from $0$ to $1$ the value of each
entry of $S(V)_{a,b}$ in position $(i',j')$, with $j+1\leq j' \leq
n'+q-1$, so that it becomes a $*1$-row. Discard such a valuation
if a row having all the entries equal to $1$ has eventually been
created, in order to maintain minimality  (remind
Lemma~\ref{lem6});

\item{$iv)$} there exists a $*1$-row $i'$, with $i+1\leq i' \leq
n'+q-1$ (see Fig.~\ref{sm40}, $(d)$).

No changes in the last $(m'+p-1)-i$ rows of $V$ allow the
insertion of the new point $p_{t+1}$.
\end{description}

\begin{figure}
\centerline{\hbox{\psfig{figure=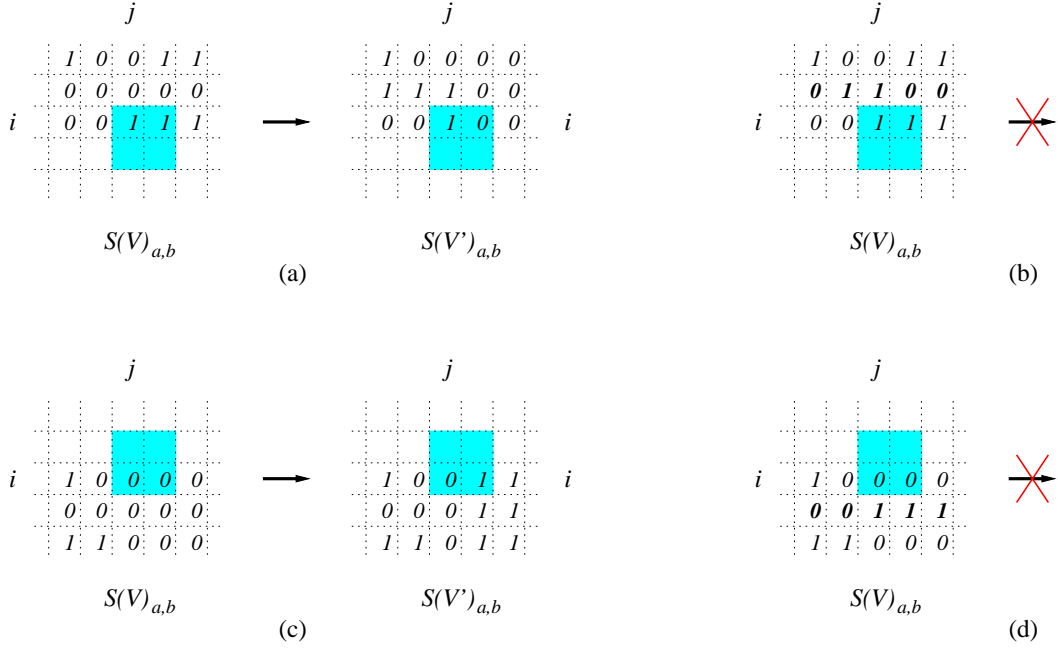,width=5.5in,clip=}}}
\caption{\label{sm40} \em \small Examples of the possible ways of
extending a valuation $V$ when adding point $p_{t+1}$. In cases
$(b)$ and $(d)$, the boldface rows prevent $V$ from being
extended.}
\end{figure}

{\em Remark:} If the point $p_{t+1}$ has value $-1$, then a
further check is needed in the analogous of case $i)$: it may
happen that there exists a $*0$-row $i'$, with $1\leq i'\leq i$,
which is turned into a $*1$-row after the addition of $p_{t+1}$,
so that a non minimal configuration is created.

The four configurations above described are exhaustive with
respect to the addition of the single point $p_{t+1}$ of value $1$
to the valuation $S(V)_{a,b}$. However, a case has not yet been
considered: it appears when two points $p_{t+1}$ and $p_{t+2}$ are
added to $S(V)_{a,b}$, under the assumption that they have
different value, and they lie in the same column.

\medskip

\item{\bf Step $t\rightarrow t+2$}: let us assume that
$p_{t+1}=(i,j)$, $p_{t+2}=(i',j)$, $S(\chi(M)_{a,b}[i,j]=1$, and
$S(\chi (M))_{a,b}[i',j]=-1$.

If it holds that:

\begin{description}

\item{$v)$} all the rows from $i$ to $i'$ are not $*1$-rows of
$S(V)_{a,b}$, and there exists a $*1$-row with index greater than
$i'$ (this last condition prevent $S(V)_{a,b}$ from being extended
to $S(V')_{a,b}$ by means of $iii)$).

We define the valuation $V'$ as follows:

for each $i\leq i_0 \leq i'$, change from $0$ to $1$ the value of
each entry of $S(V)_{a,b}$ in position $(i_0,j')$, with $j+1\leq
j' \leq n'+q-1$, so that it becomes a $*1$-row.

\end{description}
\end{description}

In the sequel, when we mention the above described cases $i)-v)$,
we intent to include also their symmetrical counterparts. It is
immediate to check that
\begin{description}

\item{-} cases $i)-v)$ extend $S(V)_{a,b}$ by adding the desired
point (or points);

\item{-} each extension of $S(V)_{a,b}$ is minimal, since no rows
or columns completely filled with entries $1$ are added;

\item{-} all the minimal valuations for the sequences
$p_1,\dots,p_{t+1}$ or $p_1,\dots,p_{t+1},p_{t+2}$ are obtained by
means of $i)-v)$.

\end{description}

\begin{figure}
\centerline{\hbox{\psfig{figure=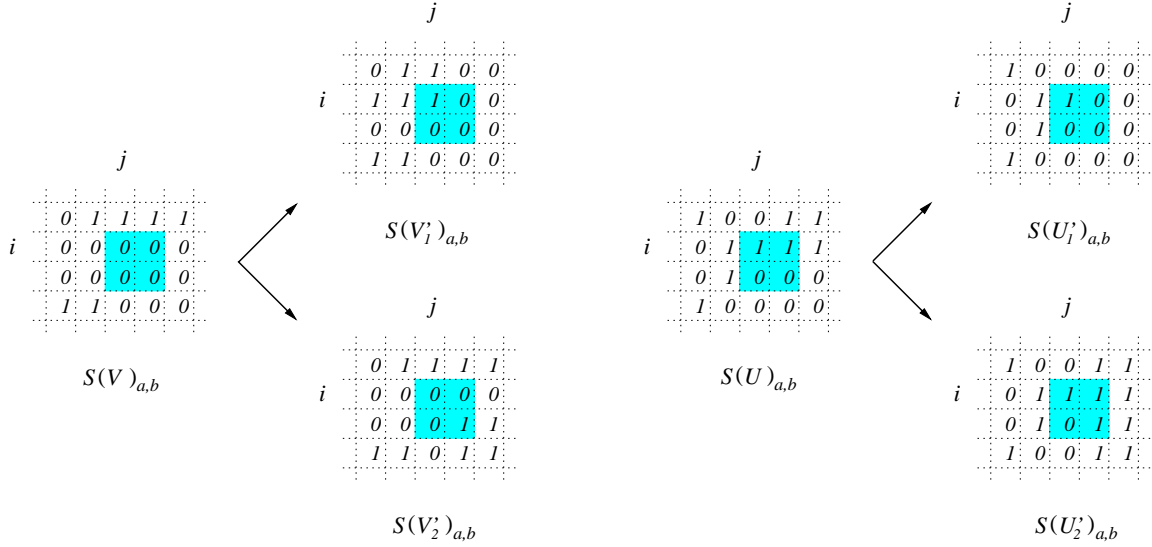,width=6in}}}
\caption{\label{sm52} \em \small The valuations $S(V)_{a,b}$
 and $S(U)_{a,b}$ extend in two different valuations. Since
 $S(U)_{a,b}$ does not contain any $0$-row, one of its extensions
(i.e. $S(U_2')_{a,b}$) is not minimal.}
\end{figure}

So, what remains to prove is that the number of the different
minimal valuations for a given matrix $S(\chi(M))_{a,b}$ is
polynomial in its dimensions $m'$ and $n'$ (and consequently in
the dimensions $m$ and $n$ of $M$). We achieve this aim by showing
that the number of valuations $S(V)_{a,b}$ which admit more than a
single minimal extension is bounded by $m'$. Some properties are
needed:

\begin{prop}\label{prop:1}
Each valuation $S(V)_{a,b}$ admits at most two different minimal
extensions both when adding a single point $p_{t+1}$ (see
Fig.\ref{sm52}), and when adding two points $p_{t+1}$ and
$p_{t+2}$, under the assumptions of $v)$.
\end{prop}

\begin{prop}\label{prop:3}
If the valuation $S(V)_{a,b}$ does not contain any $0$-row, then
it admits at most one minimal extension via $i)-v)$ (see
Fig.\ref{sm52}, valuation $U$).
\end{prop}

\begin{prop}\label{prop:4}
The $0$-rows of each valuation which extends $S(V)_{a,b}$ are a
subset of those of $S(V)_{a,b}$. Furthermore, if $S(V)_{a,b}$
extends in two minimal ways, then the two extensions do not share
any $0$-row.
\end{prop}

\begin{prop}\label{prop:5}
Two minimal valuations of $S(\chi(M))_{a,b}$ do not share any
$0$-row.
\end{prop}

The proofs of these properties directly follow from the
definitions of $i)-v)$.

Hence, Property \ref{prop:5} assures that each matrix
$S(\chi(M))_{a,b}$ has at most $m'$ different minimal valuations
containing $0$-rows. From Property \ref{prop:3}, if we add one or
two new points to them, then at most $m$ new minimal valuations
may arise. As a neat consequence, we obtain that the number of
minimal valuations of a given matrix $S(\chi(M))_{a,b}$ is
polynomial in $m'$ and $n'$, and so it is the complexity of their
reconstruction. \qed

%%%%%%%%%%%%%%%%%%%%%%%%%%  FINE PROVA %%%%%%%%%%%%%%%%%

\begin{exmp}\label{ex2}
{\small Let us find all the minimal valuations of the matrix
$S(\chi (M))_{a,b}$ depicted in Fig.~\ref{sm27}

\begin{figure}[htd]
\centerline{\hbox{\psfig{figure=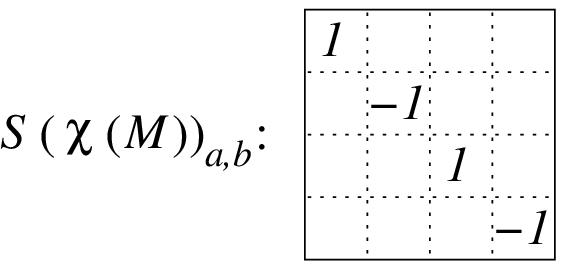,width=1.5in}}}
\caption{\label{sm27} \em \small The matrix $S(\chi(M))_{a,b}$ of
Example~\ref{ex2}. }
\end{figure}

We proceed from the leftmost entry of $S(\chi(M))_{a,b}$ different
from $0$, till the rightmost one, and we construct, step by step,
all the possible minimal valuations, as described in $i)-v)$.

The computation is represented in Fig.~\ref{sm28}, by using a tree
whose root is the matrix having all the entries equal to $0$, and
which represents the $(a,b)$ subgrid of the valuation of a
$0$-homogeneous matrix. The nodes at level $k$ are all the
possible $(a,b)$ subgrids of the valuations of the first $k$
entries different from $0$ of $S(\chi(M))_{a,b}$.

On each matrix, the highlighted cells refer to the correspondent
entry $1$ or $-1$ of $S(\chi (M))_{a,b}$ which is being
considered.

\begin{figure}[htd]
\centerline{\hbox{\psfig{figure=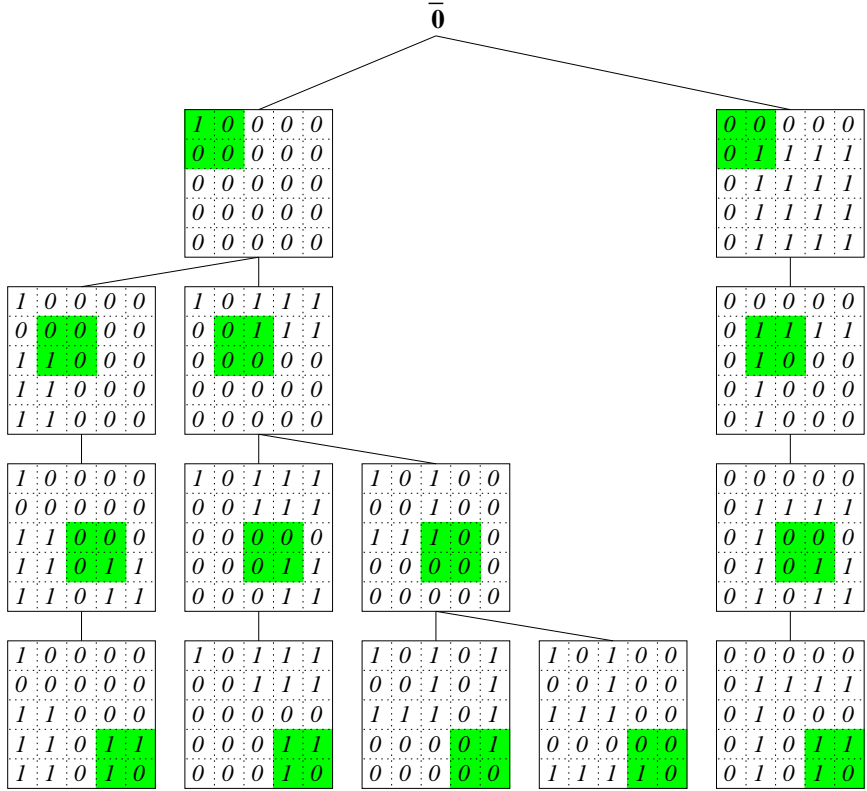,width=4.5in,clip=}}}
\caption{\label{sm28} \em \small The computation of the minimal
valuations of $S(\chi(M))_{a,b}$.}
\end{figure}

It is easy to check that any further addition of $1$ or $-1$
entries in $S(\chi(M))_{a,b}$ does not increase the number of the
minimal valuations.}

\end{exmp}

The following variant of the procedure {\sc RecSmooth} will be
used in the final reconstruction algorithm; the details of the
procedure which differ from the original ones are given:

\medskip

{\sc RecSmoothAll}$(A,p,q)$

\begin{description}

\item{\bf Input:} an integer matrix $A$ and two integers $p$ and
$q$.

\smallskip

\item{\bf Output:} a (eventually void) sequence of $m\times n$
matrices $M_1,\dots, M_k$ having elements in $\{1_P^{(1)},\dots,
1_P^{(k_{col})} ,1_Q^{(1)},\dots,1_Q^{(k_{row})},P,Q,0\}$, with
$1\leq n \leq p\times q$.

\smallskip

\item{\bf Procedure:}

\item{Step $1$:} ...

\medskip

\begin{description}
\item{Step $1.1$:} ...

\medskip

\item{Step $1.2$:} ...
\medskip

\item{Step $1.3$:} list all the possible $p\times q$ matrices such
that, for each of them

\smallskip

$i)$ the number of the entries $Q$ in its $i$-th row is
$PEnt_{\equiv p}[i]$.

\smallskip

$ii)$ the number of the entries $P$ in its $j$-th column is
$PEnt_{\equiv q}[j]$;

\smallskip

$iii)$ at least one occurrence of each entry in $1_P^{(1)},\dots,
1_P^{(k_{col})}$ and at least one occurrence of each entry in
$1_Q^{(1)},\dots , 1_Q^{(k_{row})}$ is present. Furthermore, all
the entries with pedex $P$ [resp. $Q$], and having the same index
must lie in the same column [resp. row].

{\bf If} $t\not=k$ {\bf then} set $t=t+1$, and return to
Step~$1.1$.

\end{description}

\medskip

\item{Step $2$:} let $W_1,\dots W_K$ be the output list of Step
$1$. For each $1\leq i \leq K$, use matrix $W_i$ to create the
$m\times n$ matrix $M_i$ whose entries are initialized as follows:
\begin{description}
\item{Step $2.1$:} ...

\item{Step $2.2$:} ...

\item{Step $2.3$:} for each $0\leq i' \leq p$ and $0\leq j \leq
q$,

{\bf if} $W_i[i',j]\not\in \{P,Q\}$, then set
$M_i[i',j]=W_i[i',j]$];

\item{Step $2.4$:} fill the matrix $M$ imposing the
$(0,q)$-invariance of its entries $Q,1_P^{(1)},\dots,
1_P^{(k_{col})},1_Q^{(1)},\dots,1_Q^{(k_{row})}$, and the
$(p,0)$-invariance of its entries $P,1_P^{(1)},\dots,
1_P^{(k_{col})},1_Q^{(1)},\dots,1_Q^{(k_{row})}$;
\end{description}

\medskip

\item{Step $3$:} return the sequence $M_1,\dots M_K$ as output.
\end{description}

This variant of {\sc RecSmooth} inherits its $O(m\: n)$
computational complexity.

As a final observation, one may wonder the meaning of the indexed
entries $1_P$ and $1_Q$ inside each matrix: the elements $1_P$
[resp. $1_Q$] having the same index mark the positions where a set
of entries whose rectangular scan is $1$-homogeneous, can be
placed. The choice of all the possible sets of positions for the
placement of the $k_{row}+k_{col}$-homogeneous part of $A$
constitutes a key point in the definition of the final
reconstruction algorithm which follows:

\medskip

{\sc Reconstruction}$(A,p,q)$

\begin{description}

\item{\bf Input:} an integer matrix $A$ and two integers $p$ and
$q$.

\smallskip

\item{\bf Output:} an $m\times n$ binary matrix $M$ having $A$ as
$(p,q)$ rectangular scan, if it exists, else return FAILURE.

\smallskip

\item{\bf Procedure:}

\item{Step 1:} for each $1\leq a \leq p$ and $1\leq b \leq q$,
compute the sequence of minimal valuations
%$V_{a,b}^{(1)},\dots , V_{a,b}^{(v)}$
of $S(\chi_{1,1} (A))_{a,b}$;

\smallskip

\item{Step 2:} sum in all possible ways an element from each
sequence of valuations computed in Step $1$, and let
%$M_1,\dots M_{v'}$
$M_1,\dots M_{v}$ be the obtained sequence of binary matrices;

\smallskip

\item{Step 3:} for each $1\leq t \leq v$,

\begin{description}

\item{Step $3.1$:}  compute the matrix $A_t=A-R(M_t)$;

\smallskip

\item{Step $3.2$:}  run {\sc RecSmoothAll}$(A_t,p,q)$, and let
$M_1',\dots, M_K'$ be its output. Set $t'=1$;

\smallskip

\item{Step $3.3$:}  until $t' \leq K$, compute a matrix $M$ by
merging the matrix $M_t$ and the matrix $M_{t'}'$ as follows:
initialize $M=M_t$;

for each $1\leq i \leq m$, $1\leq j \leq n$

\begin{description}

\item{\bf if} $M_{t'}'[i,j]\in \{P,Q\}$, {\bf then}

{\bf if} $M[i,j]=1$, {\bf then} set $t'=t'+1$ and return to Step
$3.3$, {\bf else} $M[i,j]=1$;

\item{\bf if} $M_{t'}'[i,j]=1_P^{(n_0)}$, with $1\leq n_0 \leq
k_{col}$, {\bf then}

{\bf if} $M[i',j]=0$, for each position $(i',j)$, with
$i'=(i)mod_p$, {\bf then} set $M[i',j]=1$, and change to $0$ the
remaining entries in column $j$ of $M_{t'}'$ having value
$1_P^{(n_0)}$ {\bf else} set $M_{t'}'[i',j]=0$, and, if no other
elements $1_P^{(n_0)}$ are in column $j$, set $t'=t'+1$ and return
to Step $3.3$;

\item{\bf if} $M_{t'}'[i,j]=1_Q^{(m_0)}$, with $1\leq m_0 \leq
k_{row}$, {\bf then}

{\bf if} $M[i,j']=0$, for each position $(i,j')$, with
$j'=(j)mod_q$, {\bf then} set $M[i,j']=1$, and change to $0$ the
remaining entries in row $i$ of $M_{t'}'$ having value
$1_Q^{(m_0)}$ {\bf else} set $M_{t'}'[i,j']=0$, and, if no other
elements $1_Q^{(m_0)}$ are in row $i$, set $t'=t'+1$ and return to
Step $3.3$;
\end{description}
{\bf Return} matrix $M$ as output;
\end{description}
\item{Step 4:} return FAILURE.
\end{description}

The correctness of the procedure is straightforward, since we
create {\em all} the possible minimal valuations for the entries
of $A$ which prevent it from being smooth, and successively, we
merge them with {\em all} the possible solutions for its remaining
smooth entries.

\begin{thm}
The problem $Reconstruction$ $(A,p,q)$ admits a solution, if and
only if the algorithm {\sc Reconstruction} $(A,p,q)$ finds it.
\end{thm}

However, one can ask whether such a search always produces an
output in an amount of time which is polynomial in the dimensions
$m$ and $n$ of the solution. The answer is given in the proof of
the following

\begin{thm}
The computational complexity of {\sc Reconstruction}$(A,p,q)$ is
polynomial in the dimension $m \times n$ of the solution.
\end{thm}

\pf The complexity of the algorithm can be computed as the sum of
the complexities of its steps, in particular:

\begin{description}

\item{Step $1$:} %for each $1\leq a \leq p$ and $1\leq b \leq q$,
%$S(\chi_{1,1}(A))_{a,b}$ can be computed in $O(m\: n)$.
Lemma~\ref{lem7} assures that the computation of all the
valuations of $S(\chi_{1,1} (A))_{a,b}$ can be performed in
polynomial time with respect to $m$ and $n$.

\item{Step $2$:} the procedure asks for summing in all possible
ways an element from each of the $p \cdot q$ sequences of
valuations created in Step $1$. Since each sum is performed in
$O(m \: n)$, then the total complexity remains polynomially
bounded by $m$ and $n$.

\item{Step $3$:} the computation of the matrix $A_t$, the
polynomial procedure {\sc RecSmoothAll} and the merging process of
$M_t$ with $M_{t'}'$ are performed a polynomial number of times,
without increasing the total complexity of the algorithm.\qed

\end{description}

\end{document}